# Reconstruction of Fractional Brownian Motion Signals From Its Sparse Samples Based on Compressive Sampling

Andriyan Bayu Suksmono

School of Electrical Engineering and Informatics
Institut Teknologi Bandung, Jl. Ganesha No10, Bandung, Indonesia
suksmono@yahoo.com; suksmono@stei.itb.ac.id

*Abstract*—This paper proposes a new fBm (fractional Brownian motion) interpolation/reconstruction method from partially known samples based on CS (Compressive Sampling). Since $1/f$ property implies power law decay of the fBm spectrum, the fBm signals should be sparse in frequency domain. This property motivates the adoption of CS in the development of the reconstruction method. Hurst parameter *H* that occurs in the power law determines the sparsity level, therefore the CS reconstruction quality of an fBm signal for a given number of known subsamples will depend on *H*. However, the proposed method does not require the information of *H* to reconstruct the fBm signal from its partial samples. The method employs DFT (Discrete Fourier Transform) as the sparsity basis and a random matrix derived from known samples positions as the projection basis. Simulated fBm signals with various values of *H* are used to show the relationship between the Hurst parameter and the reconstruction quality. Additionally, US-DJIA (*Dow Jones Industrial Average*) stock index monthly values time-series are also used to show the applicability of the proposed method to reconstruct a real-world data.

## I. INTRODUCTION

Fractional Brownian motion, or fBm [1], is a zero mean Gaussian process with statistical self-similar property. The fBm is an important model for non-stationary signals, which is well-suited to represent various kinds of natural signals, such as physiological signals, terrain surface, speech signals, internet data traffic, financial time-series, etc. The fBm was discovered by Kolmogorov and then developed and popularized by Mandelbrot. The fBm can be considered as an extension of Brownian motion into fractional dimension.

Construction of an fBm signal for an arbitrary parameters, which is known as Hurst parameter *H* or spectral parameter $\gamma = H + 1$, can be performed in space/time domain or in frequency domain. Since *H* characterizes an fBm signal, estimation of this parameter for a given signal is an important issue. The reconstruction problem will be more complicated if the fBm signal is contaminated by noise. The estimation of an fBm signal from its noisy measurement has been described in [2], [3], and [4]. In [5], an fBm equalization method and its application to improve DEM (Digital Elevation Model) reconstruction for a given InSAR (Interferometric Synthetic Aperture Radar) phase image were proposed.

Another important issue to address in the fBm reconstruction is when only a small number of data is known. In [6], a method to interpolate the fBm signal from a coarser- to a finer- grid was proposed. This kind of method can be considered as an interpolation process from regularly-spaced samples. Related to this issue, this paper extends the fBm interpolation problem further to cover non-regularly-spaced samples.

A straight-forward solution for such a reconstruction problem can be obtained by, for example, first estimating the Hurst parameter and then choose a signal, from an ensemble of *H*-parameterized fBm signals that best-fit the observation while maintaining the fractal property described by *H*. A simpler method is by linearly interpolating the blank points between known values.

In this paper, a different approach to interpolate/ reconstruct the fBm signals based on the emerging CS (Compressive Sampling) paradigm is proposed. The CS method is capable to reconstruct a signal from a few numbers of samples, when an appropriate sparsity and projection basis can be obtained. Although an fBm signal is non-sparse in space/time domain, the $1/f$ property of the fBm introduces power law decays of the Fourier coefficient, indicating that the fBm should actually a sparse signal in frequency domain. The sparsity of the fBm and CS reconstruction methods will be explored in this paper, including its potential applications.

The rest of this paper is organized as follows. Section II described the fBm properties in time/space and frequency domain. A brief review of CS theory and the proposed method will be explained in Section III. Experiments with simulated fBm signal and actual financial time series will be described in Section IV, while Section V concludes this paper.

## II. FRACTIONAL BROWNIAN MOTION SIGNAL IN TIME - AND FREQUENCY- DOMAIN

The fractional Brownian motion (fBm) is a continuous zero-mean Gaussian process that generalizes the ordinary Brownian motion. An fBm with Hurst parameter $H$, denoted by $B_H(t)$, where $0<H<1$ and initial value $B_H(0)=0$, is given by [1]

$$B_H(t) = \frac{1}{\Gamma(H+\tfrac{1}{2})}\left\{\int_{-\infty}^{0}\left(|t-s|^{H-\tfrac{1}{2}}-|s|^{H-\tfrac{1}{2}}\right)dB(s)+\int_{0}^{t}|t-s|^{H-\tfrac{1}{2}}dB(s)\right\} \quad (1)$$

Although the fBm is a non-stationary and statistically self-similar, its increments are actually stationary. The correlation function of the fBm is given by

$$E\{B_H(t)B_H(s)\} = \frac{C_H}{2}\left(|t|^{2H}+|s|^{2H}-|t-s|^{2H}\right) \quad (2)$$

where $C_H$ is an $H$-dependent coefficient given by [7]

$$C_H = \Gamma\left(1-2H\frac{\cos\pi H}{\pi H}\right) \quad (3)$$

In the frequency domain, the fBm is characterized by the well known $1/f$ profile. Denoting the frequency as $\omega$, Mandelbrot and van Ness [1] suggested that the spectral density of an fBm with Hurst parameter $H$ should proportional to $|\omega|^{-2H-1}$. Although the spectral densities of a nonstationary random function is not obviously defined, by averaging the time-dependent Wigner-Ville spectrum, Flandrin [7] showed that the average spectra over time of an fBm can be expressed as

$$S_{B_H}(\omega) \sim 1/|\omega|^{2H+1} \quad (4)$$

Various kind of signals from natural phenomena exhibit $1/f$-type spectral behaviour, among others are: geophysical-, financial-, physiological-, and biological-time series, and also the EM/electronic device-fluctuations, burst error patterns in communications channel, natural terrain texture, etc. In [4], Wornell further extended the notion of $1/f$ to cover a broader class, which is called nearly $1/f$ process.

This paper considers discrete-time fBm signals. Some methods to generate the fBm signals for a given value of $H$ have been developed, such as the midpoint random displacement or Fourier synthesis methods. Fig.1 shows realization of four different fBm signals generated by Fourier synthesis method: (a) time-domain and (b) frequency domain signals. The $H$-values of these signals, from top to bottom, are respectively 0.2, 0.4, 0.6, and 0.8. The figure shows that the signal becomes smoother as $H$ values close to 1.0, and it becomes more rough when $H$ approaching zero. Fig. 1(b) displays the localization property of the fBm spectrum; i.e., they are mostly located near the center frequency and fast decays for higher frequency components away from the center. Top to down sequence of the curves in Fig.1 (b) also shows that the spectrum are become more localized with the increasing value of $H$.

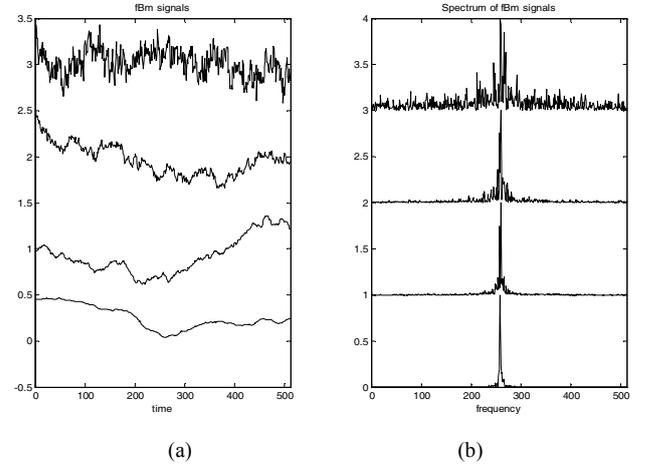

Fig.1 The fBm signals with various values of Hurst parameters: (a) time domain and (b) spectra

## III. RECONSTRUCTION OF THE FBM SIGNAL FROM ITS SPARSE SAMPLES

### A. Compressive Sampling and Signal Sparsity

Compressive Sampling/Compressed Sensing, or CS, is an emerging paradigm that unifies sampling and compression into a single process. In the conventional sampling, the Shannon theorem states that a $\Delta\omega$ bandlimited signals requires $2\times\Delta\omega$ sampling-rate for an exact reconstruction. The sampling theorem implies that the signal and sampling process is related to Fourier decomposition, where sinuoids basis is used. When different basis is employed, the requirements of the sampling rate may change. The CS extends the sampling process into more general basis. Whereas the bandwidth indicates the information content of the signal in the classical sampling, CS uses sparsity level or degree of freedom.

Underlying CS is the UP (*Uncertainty Principle*) which states that a signal cannot be simultaneously well-localized in time and frequency. The UP is closely related with the HUP (Heisenberg UP) in quantum mechanics. By expressing the HUP in term of time and energy, or time and frequency (since the quanta of energy is related to frequency by the Planck constant $h$), it can be shown that the product of the time uncertainty $\Delta t$ and its corresponding frequency uncertainty $\Delta\omega$, should be larger than a particular constant. In discrete form, the UP leads to signal recovery from highly incomplete samples [8], [9]. This principle was also extended into a more general pair of basis, instead of the time-frequency pair.

Following [10], the CS and its reconstruction are formulated as follows. Let $f$ be an $N$-length discrete time signal and $y$ be the transform of $f$ by an orthonormal (unitary) basis $\Psi=[\psi_1\ \psi_2\ ...\ \psi_N]$, i.e.,

$$f(k) = \sum_{i=1}^{N} x_i \psi_i(k) \qquad (5)$$

where $x_i = \langle f, \psi_i \rangle$ are the transform coefficients. Being a $K$-sparse signal means that by retaining only $K$-coefficients, the signal $f$ can be reconstructed with no perceptual loss. The sensing/sampling is performed incoherently by projecting the signal $f$ into an orthonormal basis $\Phi = [\phi_1\ \phi_2\ ...\ \phi_N]$, i.e.,

$$y(k) = \langle f, \phi_k \rangle \qquad (6)$$

By defining the coherence between the sensing basis $\Phi$ and sparsity basis $\Psi$ as

$$\mu(\Phi, \Psi) = \sqrt{N} \max_{1 \leq k,j \leq N} \left| \langle \phi_k, \psi_j \rangle \right|, \qquad (7)$$

then an $M$ subsamples of $f$, where [11]

$$M \geq C \cdot \mu^2(\Phi, \Psi) \cdot K \cdot \log(N) \qquad (8)$$

is sufficient (with overwhelming probability) to reconstruct $f$ using the following convex optimization

$$\min_{\tilde{x} \in R^N} \|\tilde{x}\|_{l_1} \quad subject\ to \quad y_k = \langle \phi_k, \Psi \tilde{x} \rangle,\ \forall k \in M \qquad (9)$$

The sparsity basis $\Psi$ can be any orthogonal (unitary) basis, such as wavelet, Hadamard, DCT, or DFT, based on which the signal can be represented by the least number of coefficient $K$, while the projection basis $\Phi$ should be selected to give the best coherency with $\Psi$.

The minimization $\|\tilde{x}\|_{l_1}$ means to look for the sparsest signal in the basis $\Psi$. It can also be reformulated as searching for the sparsest gradient [12] or minimization of the TV (Total Variance), which can be expressed as [9]

$$\min_{\tilde{x} \in R^N} TV(\tilde{x}) \quad subject\ to \quad y_k = \langle \phi_k, \Psi \tilde{x} \rangle,\ \forall k \in M \qquad (10)$$

In the following section, reconstruction of fBm signal from its sparse samples will be formulated as the compressive sampling problem. Some optimization tools that have been developed by CS community, such as the BP (Basis Pursuits) algorithm, *L1-magic*, *recPF*, and *twist* can be used to solve this problem.

### B. Sparse Reconstruction of The fBm Signals
#### 1) Sparsity of The fBm Signals in Frequency Domain

According to equation (4), the Fourier magnitude coefficients of an fBm signal decay as the frequency increased. The speed of the decays is determined by value of the Hurst parameter; which means that the decay will be faster for a larger value of $H$ than the smaller ones. Based on this observation, an fBm can be considered as a sparse signal in the Fourier domain, whose degree of sparsity is determined by $H$.

To quantify the sparsity, various fBm signals with four different $H$ values {0.2, 0.4, 0.6, 0.8} are generated. For each $H$ value, ten different fBm signals are generated and the percentage of dominant coefficient is counted and averaged over the ten signals. The dominant coefficient is defined as the ones whose magnitude is greater than 10% maximum. The result is displayed in Table 1.

TABLE 1
NUMBERS OF DOMINANT COEFFICIENTS IN fBM SIGNALS
FOR VARIOUS VALUES OF $H$

| Hurst Parameter ($H$) | 0.2 | 0.4 | 0.6 | 0.8 |
|---|---|---|---|---|
| Dominant Coeficients (%) | 34.3 | 3.4 | 2.2 | 1.3 |

In the Table 1, it is observed that the percentage of dominant coefficients decrease with increasing value of $H$; in particular, about 34.3% of coefficients are dominant when $H$=0.2, then it drops drastically into 3.4% when $H$=0.4, and only about 1.3% of the coefficients are dominant when $H$=0.8. The numbers in the table confirm the $H$-dependent sparsity of the fBm signals.

Expression (8) shows that the number of samples required to exactly reconstruct the signal is proportional to the sparsity $K$. In the fBm case, since the sparsity can be related to the Hurst parameter $H$, the quality of reconstructed signal for a particular value of subsample number $M$ will be better for a larger $H$ than the smaller one. Next Sections will consider two reconstruction schemes, i.e., the BP (Basis Pursuits) and the TV-minimization, to evaluate the CS-reconstruction method of sparsely sampled fBm and the relationship between $H$ and reconstruction quality.

#### 2) Sparse Reconstruction of The fBm Signals by Basis Pursuit Method

The BP (Basis Pursuit) reconstruction is a direct realization of optimization given by (9), i.e., $f$ is chosen among all possible solutions whose $L_1$ norm is minimum. To be more general, a complex-valued $N$-length discrete fBm signal $f$ will be considered, whose subsamples is denoted by $f_{sub}$. Since the fBm is sparse in frequency domain, it is reasonable to choose discrete Fourier transform as the sparsity basis $\Psi$.

In the sparse fBm reconstruction problem, the random subsamples are known. The projection basis $\Phi$ is derived from the position of known values of $f_{sub}$ by constructing a random zero-one $M \times N$ matrix whose entries is set so that $\Phi f = f_{sub}$. Then, a matrix $A$ is build by multiplying the projection basis $\Phi$ with the DFT basis $\Psi$, i.e, $A = \Phi\Psi$. Considering Fourier relationship $f = \Psi F$, where F is a vector of Fourier coefficients, the minimization is now given by

$$\min \|F\|_{l_1} \quad subject \ to \quad AF = f_{sub} \qquad (11)$$

Optimization in (11) can be solved by convex programming [13]. The CVX optimization toolbox [14] can be used to implement (11) by using the following codes:

```
cvx_begin
    variable  F(N) complex;
    minimize( norm(F,1) );
    subject to
       A*F == f_sub;
cvx_end;
```

The optimization process yield optimum solution $F^*$, which is Fourier coefficients of the fBm signal. The reconstruction of the time/space domain fBm signal is conducted by Fourier transform as follow

$$\hat{f} = \Psi F^* \qquad (12)$$

*3) Sparse Reconstruction of The fBm Signals by TV Minimization Method*

The TV minimization method has been employed to solve various problems of reconstruction from partial Fourier samples (RPFS). In [9], a phantom image simulating MRI imaging was successfully recovered from a few number of radial Fourier coefficients by using TV. Similar case has also been shown in CS-VLBI (Compressive Sensing-Very Large Base Line Interferometry) [15] for 2D problem and CS-SFCW GPR (Compressive Sensing-Stepped Frequency Continuous Wave Ground Penetrating Radar) [16] for the 1D case. The fBm reconstruction from sparse random samples can also be formulated similarly, with a little modification.

In the RFPS problems, one want to reconstruct a space/time signal $f$ from a few number of random Fourier coefficients of the signals $F_{sub}$. TV optimization assumes that $f$ is a TV-sparse or an $L_1$-smooth function. The partial Fourier operator $B$ is constructed based on the position of known subsample, which is almost similar to the previous BP case. Then, the following TV optimization is conducted

$$\min TV(f) \quad subject \ to \quad Bf = F_{sub} \qquad (13)$$

Considering $1/f$ spectrum of the fBm given by (4), the fBm can be considered as a smooth function in the transform domain. Sparse samples in spatial domain $f_{sub}$ is reconstructed by seeking the most TV-sparse Fourier transform signal $F$ as follows

$$\min TV(F) \quad subject \ to \quad B^{-1}F = f_{sub} \qquad (14)$$

where $B^{-1}$ is the inverse of partial Fourier transform. After finding the transform domain TV-sparse solution $F^*$, then the fBm signal is reconstructed by taking its inverse Fourier transform given by (12).

## IV. EXPERIMENTS AND ANALYSIS

*A. Sparse Reconstruction of Simulated fBm Signal*

In the first experiment, a 128-length complex-valued fBm signal $f$ with Hurst parameter $H$=0.75 is generated by Fourier synthesis. The signal is subsampled 4 times, yields a 32-length random subsamples signal $y$. The signal is then reconstructed using the CS-BP and CS-TV, where CVX and $L_1$-Magic respectively, are employed.

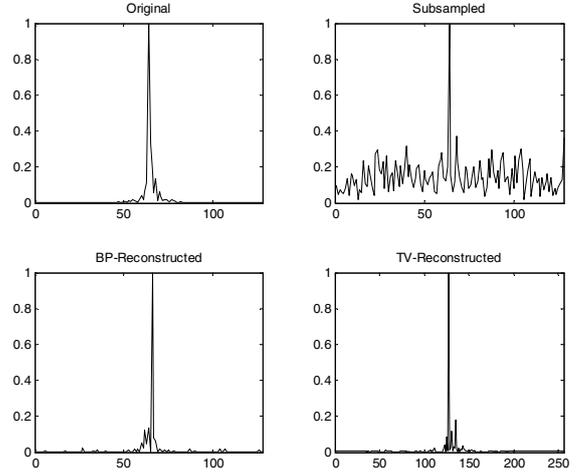

Fig. 2 Spectra of the original, subsampled, and CS-reconstructed fBm signals.

Comparisons of signal spectra are displayed in Figure 2. The spectrum of the subsamples given in the right upper part of Fig.2 shows that random subsampling makes the originally smooth and sparse signal in the left-upper part spread over the frequency shown in the right-upper part. The CS reconstruction seek for a solution with the sparsest (least $L_1$ coefficient or smoothest spectrum) that consistent with observation given by the subsamples. The solution given in the lower part of the figure, i.e. BP in the left and TV in the right parts, show the desired spectra of the solution. After obtaining the most suitable specrum/transform coefficients, the signal can be reconstructed whose results is displayed in Figure 3.

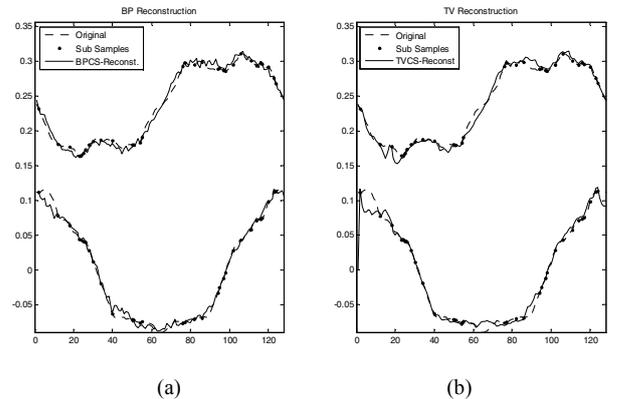

(a) (b)
Fig. 3 Reconstruction of complex-valued fBm signal using CS. The upper curves are imaginary parts, while lower ones are the real parts: (a) using BP and (b) using TV minimization

In the Fig. 3, the upper curves show imaginary parts, while the lower ones are real parts, with (a) shows BP reconstruction and (b) displays TV reconstruction result. The circle in the curve show random subsamples while the dotted lines are original fBm signals. Reconstruction results depicted by the solid curves well-fitted the original signals. Signal to noise ratio (SNR) measurements indicate that the BP yields about 19.7 dB while TV achieves 15.5 dB. The curves in the figure demonstrate that CS method successfully reconstruct the fBm signal from its sparse samples.

TABLE 2
RECONSTRUCTION QUALITY FOR VARIOUS VALUE OF $H$

| No | Sub sampling Factor | Hurst Parameter $H$ | SNR CS-BP (dB) | SNR CS-TV (dB) |
|---|---|---|---|---|
| 1 | 2 | 0.2 | 2.69 | 3.79 |
| 2 | | 0.4 | 8.29 | 8.42 |
| 3 | | 0.6 | 14.54 | 13.80 |
| 4 | | 0.8 | 24.85 | 19.19 |
| 5 | 4 | 0.2 | -0.23 | 0.12 |
| 6 | | 0.4 | 4.57 | 4.65 |
| 7 | | 0.6 | 10.99 | 11.28 |
| 8 | | 0.8 | 15.64 | 14.24 |

The effect of Hurst parameter in the reconstruction quality is displayed in Table 2. The performance is measured using SNR and two subsampling factors are considered. The subsampling factor is defined as the ratio of original number of samples with the number of subsamples. Hurst parameter is varied from 0.2 to 0.8, representing a more localized and smooth spectrum, as discussed in the previous section, with increasing $H$. Reconstruction is conducted using both the BP and TV minimization, indicated as CS-BP and CS-TV subsequently in the table. The SNR performances listed in the last two columns, which are average of 10 trial for each values, consistent with the prediction that for a given number of subsamples, the reconstruction quality will become better with the increasing value of $H$.

### B. Sparse Reconstruction of A Financial Time Seris Data

Financial time series are nonstationary and have self-similar property like the fBm signals. In this section, the capability of sparse fBm reconstruction method to an actual financial time series data is demonstrated. A 512 length US-DJIA (United States - Dow Jones Industrial Average) stock index monthly values, started from January 1945 up to August 1987, is used. This sequence is subsampled 4 times, giving 128 length random samples. Then, both of the CS-BP and TV-BP algorithm are employed to recover the full-length signal.

Figure 4 displays the spectra of the US-DJIA. The top-left part shows the original spectrum, displaying a smoother one compared to fBm spectrum of $H$=0.75 displayed previously in Fig.2. Subsampling yields the spreading of the spectrum displayed in the top-right part. After recovery by the CS methods, the spectrum changes into ones showed in the bottom parts, CS-BP is given in the left while CS-TV is located in the right part.

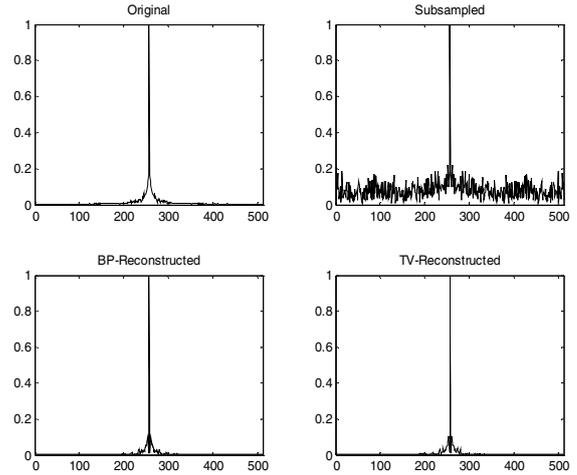

Fig. 4 Spectra of 512-length DJIA samples, its subsamples, and recovered spectrum by CS method

Figure 5 shows the reconstruction result of the US-DJIA: (a) using CS-BP and (b) using CS-TV. Measurement of the average SNR over ten times experiments shows that both of these methods are comparable, CS-BP give 23.2 dB while the CS-TV achieved 22.4 dB. They even perform better than the ones in the simulation. Observation of the spectra in Fig.4 shows may explain the reason, the US-DJIA spectra is more localized than simulated fBm signals. It is also found that for a given $H$ values assumed constant along the signal, reconstruction of longer signals gives better SNR than the shorter ones. The detail analysis and experiments regarding this issue will be explained in the future research.

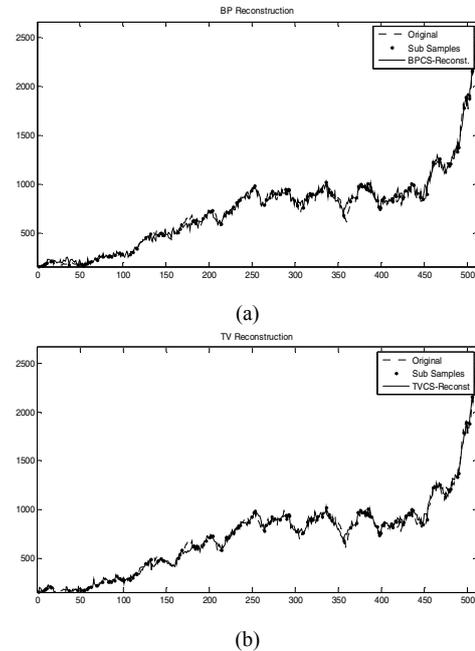

Fig. 5 Reconstruction result of the *US-DJIA* stock index monthly values from its random subsamples

## V. Conclusions

The fBm signal interpolation/reconstruction from its subsamples has been formulated as a compressive sampling problem. Two reconstructions approached have been presented, the CS-BP and CS- TV. The CS-BP seeks for the sparsest transform coefficient, whereas the CS-TV looks for the smoothest one. Their performances are comparable for simulated fBm signals as well as real-life financial time series. Frequency domain representation of fBm, indicating *H*-dependent power decay of transform coefficients, implies that for a given number of subsamples, reconstruction of large *H* will give a better result than the lower ones. The capability of CS in reconstructing real-life financial time-series, i.e. the US-DJIA, demonstrates the applicability of the proposed method.

## Acknowledgements

This work is supported by ITB PRI (*Penguatan Riset Institusi*) Research Grant 2010.